\documentclass[12pt,a4paper]{article}
\usepackage[utf8]{inputenc}
\usepackage{hyperref}
\usepackage{amsmath}
\usepackage{amsfonts}
\usepackage{amssymb}
\usepackage{mathtools}
\usepackage[margin=10pt,font=small,labelfont=bf]{caption}
\usepackage{graphicx}
\usepackage{algorithm}
\usepackage{booktabs}
\usepackage{algorithmicx}
\usepackage{algpseudocode}
\usepackage[a4paper, total={6in, 8in}]{geometry}
\usepackage{color}

\graphicspath{{../plots/}}
\usepackage{authblk}

\newcommand{\R}{\mathbb{R}}
\newcommand{\N}{\mathbb{N}}
\newcommand{\F}{\mathcal{F}}
\newcommand\Tstrut{\rule{0pt}{2.6ex}}         
\newcommand\Bstrut{\rule[-0.9ex]{0pt}{0pt}}   

\providecommand{\keywords}[1]
{
  \small	
  \textbf{\textit{Keywords---}} #1
}

\title{An exponential integrator/WENO discretization for sonic-boom simulation on modern computer hardware}
\author[1]{L. Einkemmer}
\author[1]{ A. Ostermann}
\author[1]{M. Residori}
\affil[1]{Department of Mathematics, University of Innsbruck, Austria}
\begin{document}
\maketitle

\begin{abstract}
Recently a splitting approach has been presented for the simulation of sonic-boom propagation. Splitting methods allow one to divide complicated partial differential equations into simpler parts that are solved by specifically tailored numerical schemes.
The present work proposes a second order exponential integrator for the numerical solution of sonic-boom propagation modelled through a dispersive equation with Burgers' nonlinearity. The linear terms are efficiently solved in frequency space through FFT, while the nonlinear terms are efficiently solved by a WENO scheme. The numerical method is designed to be highly parallelisable and therefore takes full advantage of modern computer hardware.
The new approach also improves the accuracy compared to the splitting method and it reduces oscillations. The enclosed numerical results illustrate that parallelisation on a CPU results in a speedup of 22 times faster than the straightforward sequential version. The GPU implementation further accelerates the runtime by a factor 3, which improves to 5 when single precision is used instead of double precision. 
\end{abstract}


\section{Introduction}
\label{sec:intro}

Sonic-booms are acoustic waves generated by supersonic planes when they fly faster than the speed of sound. Sonic-booms are heard as two loud bangs that are close together. These bangs not only annoy the population, but they can also potentially damage building facades. For these reasons supersonic planes are limited to military use and commercial supersonic flights are still not possible. Starting from the 60s many studies have been conducted in order to design the shape of the aircraft so that the generation of sonic-booms is minimized or eliminated; we refer the reader to~\cite{Alonso12,howe05,Boccadoro12, plotkin89, seebass72}. The theory finds partial confirmation in physical experiments~\cite{lipkens02}. However, conducting real experiments turns out to be extremely expensive, as the shape of the aircraft cannot be changed easily. Therefore, the need arises to perform numerical simulations that aim to model the propagation of acoustic waves generated by supersonic planes. 

It is well known that the airflow over a supersonic aircraft generates a pressure disturbance. The acoustic wave originating from the pressure disturbance will evolve in a $N$ shaped wave when it propagates ``far'' enough from the aircraft. Typically, $N$-waves appear after that acoustic waves generated by the pressure disturbance have propagated ten body lengths away from the aircraft. This distance is usually referred to as the mid field. The $N$-wave is a mathematical model to the loud bangs perceived by humans. A first issue is then to solve an inverse problem in order to predict the formation of $N$-waves in the mid field starting by a pressure disturbance generated by the aircraft geometry in the near field (i.e. the acoustic wave generated shortly after the plane passed), see~\cite{Zuazua16, Alonso12}. A second problem is to study the propagation of the resulting $N$-wave from the mid field to the ground (far field). In this paper we consider only the second problem, which is of key importance in any optimization algorithm as the final goal is to model the shape of the aircraft in such a way that $N$-waves do not appear or are mitigated by the time that sonic-booms reach the ground. This turns into modelling the sonic-boom propagation through a partial differential equation and solving it several times. We focus our attention on the evolution of acoustic waves from the mid field to the far field. To do so, we take as initial value the $N$-wave and we simulate its propagation into the far field.

In our study, we follow the mathematical model proposed in~\cite{blanc11} in order to simulate nonlinear and diffraction effects in the propagation of $N$-waves. This consists in solving numerically  a Khokhlov-Zabolotskaya-Kuznetsov (KZK)-type equation, which is done in~\cite{blanc11} by splitting methods. In this paper we introduce a new approach based on exponential integrators to solve efficiently the aforementioned partial differential equation. Exponential integrators have been successfully employed to solve various partial differential equations; we refer the interested reader to~\cite{einkemmer20, einkemmer18,ostermann05,ostermann10}. The main idea is to solve the linear part of the KZK-equation in frequency space with the help of the fast Fourier transform (FFT) and to discretise the nonlinear terms with a weighted essentially non-oscillatory (WENO) scheme. The latter scheme is very well known in the literature and is widely used in computational fluid dynamics for the numerical solution of hyperbolic conservation laws, see~\cite{ shu96,shu98,shu99}. The new approach brings multiple advantages such as the reduction of oscillations, less number of operations (asymptotically speaking) and acceleration due to parallelisation. 

The goal of this work is to fully exploit modern computer hardware such as GPUs to drastically reduce the computational cost of the numerical simulations. There is a flourishing literature about use of GPUs in order to accelerate scientific computations, e.g.~\cite{Burau10,einkemmer20a,gainullin15,griebel10,einkemmer19}. Reducing the simulation time is of great importance when we aim to model physical phenomena. Efficient implementations allow us to consider, for example, more grid points to obtain numerical solutions that are closer to reality. To provide numerical algorithms that take full advantage of parallel architectures is therefore of great practical interest. To achieve good performance it is very important to design the numerical scheme keeping the parallel paradigm in mind. Indeed, a sequential algorithm will probably run much faster on a single, powerful, single-core CPU rather than on a GPU. This is due to the fact that the number of operations per seconds on each core is much lower on a GPU compared to a CPU, see~\cite{feng19}. Increase in performance is only possible when the parallel architecture of a GPU is fully used. As mentioned above, the approach proposed in this paper is based on exponential integrators. Examples of their implementation on GPUs can be found e.g. in~\cite{kandolf18,caliari20,einkemmer13,shanks20}.

The paper is structured as follows: in section~\ref{sec:mat} we describe the mathematical model of $N$-waves propagation in randomly inhomogeneous media. In section~\ref{sec:na} we briefly recall the approach based on splitting methods and introduce a new algorithm based on exponential integrators. The two methods are discussed and compared.  In section~\ref{sec:ne} we present some numerical experiments that support the theoretical derivations and illustrate the performance of the proposed numerical scheme on parallel architectures.

\section{Mathematical model}
\label{sec:mat}
In this work we consider a mathematical model for the description of sound propagation. The model is based on a nonlinear dispersive partial differential equation that takes into account effects of the turbulent velocity field and describes the evolution of the acoustic pressure. For more details we refer the reader to~\cite{blanc11,khokhlova06}. Henceforth, we consider the following partial differential equation in the unknown  $V = V(\sigma,\rho,\theta)$:
\begin{multline}
\label{eq0}
\partial_{\theta}\bigg(\partial_{\sigma}V(\sigma,\rho,\theta)  - \frac{B}{2}\,\partial_{\theta}V^2(\sigma,\rho,\theta)  - A\,\partial^2_{\theta} V(\sigma,\rho,\theta) \\
- 2\pi U_{\parallel}(\sigma,\rho) \partial_{\theta} V(\sigma,\rho,\theta)   + U_{\perp}(\sigma,\rho)\partial_{\rho} V(\sigma,\rho,\theta)\bigg) = \frac{1}{4\pi} \partial^2_{\rho} V(\sigma,\rho,\theta)
\end{multline}
with $0<B$, $0\leq A\in\R$ constants and $U_{\parallel} = U_{\parallel}(\sigma,\rho)$, $U_{\perp} = U_{\perp}(\sigma,\rho)$ variable coefficients. Eq.~\eqref{eq0} is a KZK-type equation in dimensionless form that models the acoustic wave propagation in inhomogeneous medium. The unknown $V$ is the acoustic pressure normalized with respect to the initial pulse amplitude. The variables $\sigma$ and $\rho$ are the propagation distance and the transverse coordinate, both normalized with respect to the initial pulse length. The variable $\theta$ is the time normalized with respect to the initial pulse duration. Further information about the model and parameters are given in~\cite{blanc11}.

The variable coefficients $U_{\parallel}$ and $U_{\perp}$ are the first and second component of a two dimensional isotropic random velocity field
 \[
U = \frac{1}{c_0}\begin{bmatrix} U_{\parallel}\\ U_{\perp}\end{bmatrix},
\]
where $c_0\sim 343\,$m/s denotes the ambient sound speed.
The velocity field $U$ is computed by following the approach given in~\cite{blanc11}. At a given point $r = \lambda (\sigma,\rho)$ the velocity field is given by the sum of $N$ random modes through the formula
\begin{align*}
& U(r) = \sum_{n=1}^{N} \widetilde{U}(K_n)\cos(K_n\cdot r + \phi_n), \\
& \widetilde{U}(K_n)\cdot K_n = 0,
\end{align*}
where  ``$\cdot$'' denotes  the scalar product in $\mathbb{R}^2$. 
The angle $\phi_n$ is the phase of the $n$th  mode and $K_n$ is the wave vector given by
\[
K_n = \lvert K_n\rvert \begin{bmatrix}
\cos \theta_n\\
\sin \theta_n
\end{bmatrix},
\]
where $\theta_n$ is the angle between $K_n$ and the $\sigma$-axis. Both $\phi_n$ and $\theta_n$ are elements of two independent random sequences uniformly distributed in $[0,2\pi]$. The wavenumbers $\lvert K_n\rvert$ are equispaced in an interval $[K_{\min}, K_{\max}]$. The amplitude $\lvert \widetilde{U}(K_n)\rvert$ is related to the Gaussian energy spectrum 
\[
E(K) = \frac{1}{8}\sigma_u^2K^3L^4 \exp\left(-\left(\frac{KL}{2}\right)^2\right)
\]
through the formula
\[
\lvert \widetilde{U}(K_n)\rvert = \sqrt{\frac{E(\lvert K_n\rvert)}{N}}.
\] 
In this work we set $[K_{\min},K_{\max}]  = [0.1/L, 9.0/L]$, where $L=4\lambda$ is the length scale.  
The parameter $\sigma_{u}$ is set to $3\;$m/s and $\lambda=T_0\, c_0 $, where $T_0=2\cdot 10^{-2}\;$s is the initial pulse duration. A similar setting is adopted in~\cite{blanc11}. This gives fluctuations of the variable coefficients $U_{\parallel}$, $U_{\perp}$ so that 
\[
\lVert U_{\parallel}\rVert_{\infty}, \lVert U_{\perp}\rVert_{\infty} \leq 0.05.
\]
The information will be later useful to estimate the CFL conditions of the proposed numerical schemes. A pseudo-code for the generation of the inhomogeneous velocity fields $U_{\parallel}$ and $U_{\perp}$ is given in Algorithm~\ref{alg1}. 

\begin{algorithm}
\caption{Inhomogeneous velocity fields generator}
\label{alg1}
\begin{algorithmic}
\small
\State Generate $\{\phi_n\}$, $\{\theta_n\}$ random sequences uniformly distributed in $[0,2\pi]$ for $1\leq n\leq N$;
\State Generate $\{K_n\}$ sequence of equidistant wavenumbers in $[K_{\min},K_{\max}]$ for $1\leq n\leq N$;
\State Compute $\lvert \widetilde{U}(K_n) \rvert: = \sqrt{\frac{E(\lvert K_n\rvert)}{N}},\quad 1\leq n\leq N$;
\State  Construct $\widetilde{U}(K_n) = \begin{bmatrix}\widetilde{U}_{n,1} \\ \widetilde{U}_{n,2}\end{bmatrix} = \lvert \widetilde{U}(K_n)\rvert \begin{bmatrix} -\sin\theta_n \\ \cos\theta_n\end{bmatrix}$;
\State  $U_{\parallel}(\sigma_i,\rho_j) = \sum_n\left(\widetilde{U}_{n,1} \cos\left(K_n\lambda \sqrt{\sigma_i^2+\rho_j^2}\cos(\theta_n-\mathrm{arctan2}(\rho_j,\sigma_i))+\phi_n\right)\right)$;
\State  $U_{\perp}(\sigma_i,\rho_j) = \sum_n\left(\widetilde{U}_{n,2}\cos\left(K_n\lambda \sqrt{\sigma_i^2+\rho_j^2}\cos(\theta_n-\mathrm{arctan2}(\rho_j,\sigma_i))+\phi_n\right)\right) $;\vskip 3mm\hrule\vskip 2mm 
\State \scriptsize{The function $\mathrm{arctan2}(\rho,\sigma)$ returns the angle of the vector $(\sigma,\rho)$ with respect to the $\sigma$-axis.}
\end{algorithmic}
\end{algorithm}

\section{Numerical approach}
In this section we describe two different numerical approaches for the solution of~\eqref{eq0}. For the sake of comparison we present first a splitting approach following the one given in~\cite{blanc11}. Then, we device a new approach based on exponential integrators and WENO schemes. This new approach requires (asymptotically) a smaller number of machine operations. Moreover, the method is of second order in the variable $\sigma$. This improves the convergence rate with respect to the splitting approach, which was of first order only. Finally, we observe numerically that the second approach has smaller oscillations (in amplitude) compared to the splitting approach for long propagation distances. In the following the two methods are described in detail and compared in terms of computational cost.
\label{sec:na}
\subsection{Splitting method}
\label{subsec:spl}
A possible numerical approach is given by the Lie-Trotter splitting. This method consists in dividing~\eqref{eq0} in sub-problems each of them modelling a single physical effect. 

Before we proceed to describe the numerical scheme, we transform~\eqref{eq0} in order to obtain an evolution equation in the variable $\sigma$. To do so, we simply integrate both sides of~\eqref{eq0} from $\theta_{\min}$ to $\theta$. Therefore, we obtain 
\begin{multline}
\label{eq1}
\partial_{\sigma}V(\sigma,\rho,\theta) = \frac{1}{4\pi}  \int_{\theta_{\min}}^{\theta}\partial^2_{\rho} V(\sigma,\rho,\tilde{\theta})\,\mathrm{d}\tilde{\theta} + \frac{B}{2}\,\partial_{\theta}V^2(\sigma,\rho,\theta) \\
+ A\,\partial^2_{\theta} V(\sigma,\rho,\theta) + 2\pi U_{\parallel}(\sigma,\rho) \partial_{\theta} V(\sigma,\rho,\theta) - U_{\perp}(\sigma,\rho)\partial_{\rho} V(\sigma,\rho,\theta).
\end{multline}
In~\eqref{eq1} we assumed $\partial_{\sigma}V(\sigma,\rho,\theta_{\min})=0$ for every $\sigma$ and every $\rho$. This assumption holds true if the domain is chosen large enough with respect to $\theta$. Indeed, in this case the initial data do not evolve at the boundaries (or the effects at the boundaries are negligible). Therefore, boundary conditions do not play a significant role in the numerical simulations. We assume homogeneous Neumann boundary conditions for $\rho$ and periodic boundary conditions for $\theta$, in the same spirit as in~\cite{blanc11}.

We set $\sigma$ as the ``marching'' direction, also known as artificial time. Let $0\leq \sigma\leq \Sigma$, $N_{\sigma}\in\N$, $\Delta\sigma= \Sigma/N_{\sigma}$ and $\sigma^n = n\Delta\sigma$, $n = 0,1,\dots ,N_{\sigma}$ be the uniform discretization of the variable $\sigma$. We split up~\eqref{eq0} into five equations as follows:
\begin{alignat}{2}
\partial_{\sigma}V &= \frac{1}{4\pi}  \int_{\theta_{\min}}^{\theta} \partial^2_{\rho} V(\sigma,\rho,\tilde{\theta})\,\mathrm{d}\tilde{\theta}\qquad\quad && \text{\small{(Diffraction)}}  \label{eq1a}\\[10pt]
\partial_{\sigma}V &=   \frac{B}{2}\;\partial_{\theta}V^2 && \text{\small{(Nonlinearity)}} \label{eq1b}\\[10pt]
\partial_{\sigma}V &=  2\pi U_{\parallel}(\sigma,\rho)\; \partial_{\theta} V  && \text{\small{(Axial convection)}}\label{eq1c}\\[10pt]
\partial_{\sigma}V &=  A\;\partial^2_{\theta} V && \text{\small{(Absorption)}} \label{eq1d}\\[10pt]
\partial_{\sigma}V &=   - U_{\perp}(\sigma,\rho)\;\partial_{\rho} V  && \text{\small{(Transverse convection)}}\label{eq1e}
\end{alignat}
Then, starting from an approximation $V^n(\rho,\theta)$ to the the solution of~\eqref{eq1} at $\sigma = \sigma^n$, the solution $V(\sigma,\rho,\theta)$ at $\sigma=\sigma^n+\Delta\sigma$ is approximated by
\[
V^{n+1}(\rho,\theta) = \varphi_{\Delta\sigma}^{[5]}\circ \varphi_{\Delta\sigma}^{[4]} \circ \varphi_{\Delta\sigma}^{[3]} \circ \varphi_{\Delta\sigma}^{[2]} \circ \varphi_{\Delta\sigma}^{[1]}\, V^n(\rho,\theta),
\]
where $\varphi^{[i]}_{\Delta\sigma}$, $i=1,\dots,5$ are the flows of the initial value problems associated to \eqref{eq1a}--\eqref{eq1e}, respectively. The solution of each sub-problem is approximated by different numerical schemes that are tailored to the considered sub-problem. We remark that the Lie-Trotter splitting is a method of first order in $\sigma$. This might be insufficient for certain applications. However, if needed the scheme could be generalized to second order which increases the computational cost by approximately a factor of two.

\subsection{Full discretization of single flows}
\label{fdspl}
The numerical schemes for the single sub-problems are described in the following.  We adopt a uniform discretization both in the variables $\rho_{\min}\leq \rho\leq\rho_{\max}$ and $\theta_{\min}\leq\theta\leq\theta_{\max}$ and denote by $V^n_{j,k}$ the numerical approximation of $V(\sigma^n,\rho_j,\theta_k)$. Let 
\[
\Delta\rho = \frac{\rho_{\max}-\rho_{\min}}{N_{\rho}}\quad \text{and}\quad \Delta\theta = \frac{\theta_{\max}-\theta_{\min}}{N_{\theta}},
\quad N_{\rho},N_{\theta}\in\N,
\]
be the grid sizes for $\rho$, $\theta$, respectively. Then, the sub-problems are discretized as follows. 
\vspace{3mm}

\textbf{Diffraction}. Equation~\eqref{eq1a} is a diffraction equation and will be solved by a Crank--Nicolson finite difference scheme combined with the trapezoidal rule in $\theta$:
\begin{equation}
\label{eq12}
\frac{V^{n+1}_{j,k}-V^n_{j,k}}{\Delta\sigma} = \frac{1}{4\pi} \frac{\Delta\theta}{2}\sideset{}{^*}\sum_{l=0}^k \left(\frac{V^{n}_{j+1,l}- 2 V^n_{j,l} + V^n_{j-1,l}}{\Delta\rho^2} + \frac{V^{n+1}_{j+1,l}-2 V^{n+1}_{j,l} + V^{n+1}_{j-1,l}}{\Delta\rho^2}\right),
\end{equation}
where 
\[
\sideset{}{^*}\sum_{l=0}^k u_l = \frac{u_0}{2} + \sum_{l=1}^{k-1} u_l + \frac{u_k}{2}.
\]

\textbf{Nonlinearity}. Equation~\eqref{eq1b} is a Burgers' equation, responsible for the nonlinear effects. We employ the Godunov method, see~\cite{leVeque92}, which is conservative. The discretization is given by
\begin{equation}
\frac{V^{n+1}_{j,k}-V^n_{j,k}}{\Delta\sigma} = -\frac{F\left(V^n_{j,k},V^n_{j,k+1}\right)-F\left(V^n_{j,k-1},V^n_{j,k}\right)}{\Delta\theta},
\end{equation}
with 
\begin{equation*}
F(u_l,u_r) = \begin{cases}
\min_{u_l\leq u\leq u_r} \left(-{\frac{B}{2} u^2}\right),\quad \text{if } u_l\leq u_r,\\
\max_{u_r\leq u\leq u_l} \left(-{\frac{B}{2} u^2}\right),\quad \text{if } u_l > u_r.
\end{cases}
\end{equation*}

\textbf{Axial convection} and \textbf{Absorption}. Equations~\eqref{eq1c} and~\eqref{eq1d} model the axial convection and acoustic absorption, respectively. The solutions are computed in frequency space. Let us represent $V$ by its Fourier series in the variable $\theta$:
\[
V(\sigma,\rho,\theta) = \sum_{m} \hat{v}_m(\sigma,\rho) \exp\left(\mathrm{i}2\pi m\,\frac{\theta-\theta_{\min}}{\theta_{\max}-\theta_{\min}}\right).
\]
Then, the solution at $\sigma+\Delta \sigma$ is given by 
\begin{equation}
V(\sigma+\Delta\sigma,\rho,\theta) = \sum_m \hat{v}_m(\sigma+\Delta\sigma,\rho) \exp\left(\mathrm{i}2\pi m\frac{\theta-\theta_{\min}}{\theta_{\max}-\theta_{\min}}\right),
\end{equation}
where 
\begin{equation}
\hat{v}_m(\sigma+\Delta\sigma,\rho)  = \hat{v}_m(\sigma,\rho)\,\exp\left(\mathrm{i}\frac{2\pi m}{\theta_{\max}-\theta_{\min}} 2\pi \int_{\sigma}^{\sigma+\Delta\sigma} \!\!\!U_{\parallel}(\tilde{\sigma},\rho)\,\mathrm{d}\tilde{\sigma}-A \left(\frac{2\pi m}{\theta_{\max}-\theta_{\min}}\right)^2\!\! \Delta\sigma \right).
\end{equation}

\textbf{Transverse convection}. The last equation~\eqref{eq1e} is the transverse convection. The solution is obtained by a Lax--Wendroff method, which is conservative and of second order both in $\sigma$ and $\rho$. The numerical scheme is obtain as follows. We compute a Taylor expansion of $V$ in the variable $\sigma$:
\[
V(\sigma+\Delta\sigma,\rho) = V(\sigma,\rho) + \partial_{\rho}V(\sigma,\rho) \Delta\sigma +\partial_{\sigma}^2V(\sigma,\rho)\frac{\Delta\sigma^2}{2} + \mathcal{O}(\Delta\sigma^3)
\]
and note that 
\[
 \begin{split}
\partial_{\sigma} V = -U_{\perp}\partial_{\rho} V\quad\text{and}\quad \partial_{\sigma}^2 V = -\partial_{\sigma}(U_{\perp} \partial_{\rho}V) &= -\partial_{\sigma}U_{\perp} \partial_{\rho}V - U_{\perp}\partial_{\sigma\rho} V\\
& = -\partial_{\sigma}U_{\perp} \partial_{\rho}V + U_{\perp}\partial_{\rho} U_{\perp}\partial_{\sigma} V + U_{\perp}^2\partial^2_{\rho} V.
\end{split}
\]
Inserting $\partial_{\sigma} V$ and $\partial_{\sigma}^2 V$ in the Taylor expansion and approximating the derivatives in $\rho$ by centred finite differences gives the Lax--Wendroff scheme:
\begin{equation}
\begin{split}
 \frac{V^{n+1}_{j,k}-V^n_{j,k}}{\Delta\sigma} = &- \left(U_{\perp}(\sigma^n,\rho_j)\, \frac{V^n_{j+1,k} -V^n_{j-1,k}}{2\Delta\rho}\right)  \\
& -\frac{\Delta\sigma}{2}  \left(\partial_{\sigma}U_{\perp}(\sigma,\rho_j)\big|_{\sigma=\sigma^n} \frac{V^n_{j+1,k}-V^n_{j-1,k}}{2\Delta\rho}\right)  \\
 & + \frac{\Delta\sigma}{2}  \left(U_{\perp}(\sigma^n,\rho_j)\, \partial_{\rho} U_{\perp}(\sigma^n,\rho)\big|_{\rho=\rho_j} \frac{V^n_{j+1,k}-V^n_{j-1,k}}{2\Delta\rho}\right)\\
& +\frac{\Delta\sigma}{2}  \left(U_{\perp}^2(\sigma^n,\rho_j)\, \frac{V_{j+1,k}^n -2V^n_{j,k}+V^n_{j-1,k}}{\Delta\rho^2} \right) .       
\end{split}
\end{equation}

The global numerical scheme is of first order in $\sigma$ and of second order in $(\rho,\theta)$. We remark that~\eqref{eq1b} and~\eqref{eq1e} are solved by explicit conservative schemes. Therefore, a CFL condition has to be satisfied. In particular,~\eqref{eq1b} and~\eqref{eq1e} give 
\begin{equation}
\label{eq6}
 B\max_n \lVert V^n\rVert_{L^{\infty}(\rho,\theta)}\,\Delta\sigma \leq \Delta\theta\quad\text{and}\quad \lVert U_{\perp}\rVert_{L^{\infty}(\sigma,\rho)}\, \Delta\sigma \leq \Delta\rho,
\end{equation}
respectively.
These CFL conditions are not too restrictive, indeed, we have 
\[
B\leq 5\cdot10^{-2},\quad \lVert U_{\perp}\rVert_{L^{\infty}}\leq 5\cdot 10^{-2}\quad \text{and}\quad \max_n \lVert V^n(\rho,\theta)\rVert_{L^{\infty}} \leq 5.
\]
The second bound is obtained numerically.  The absorption coefficient $A$ in sonic-boom simulation is typically of size $7\cdot 10^{-6}$. For the above numerical scheme a stronger absorption coefficient is needed in order to ensure stability of the solution. Therefore, we set $A=3.4\cdot 10^{-4}$. 

The algorithm proceeds sequentially to solve~\eqref{eq1a}--\eqref{eq1e}. In the following we discuss the computational cost of the single steps.
\begin{itemize}
\item Step~\eqref{eq1a} is the most expensive one and requires $\mathcal{O}\left(N_{\rho}N_{\theta}^2\right)$ operations. 
\item The steps~\eqref{eq1b} and~\eqref{eq1e} are solved in $\mathcal{O}(N_{\rho}N_{\theta}$) operations. \item Steps \eqref{eq1c} and~\eqref{eq1d} are solved in $\mathcal{O}(N_{\rho}N_{\theta}\log N_{\theta})$ operations. 
\end{itemize}

Advancing the numerical solution from $V^n$ to $V^{n+1}$ requires  $\mathcal{O}\left(N_{\rho}N_{\theta}^2\right)$ operations. One of the main disadvantages of the proposed splitting approach is that it is not suitable for parallelisation. The main obstacle is given by the numerical scheme~\eqref{eq12}. Indeed, this numerical scheme approximates the integral by a sum which is a non-local operator. This means that the numerical solution at stage $k$ cannot be computed before all the solutions till stage $k-1$ are computed. Therefore, parallelising this process is not possible.  Further, the absorption parameter $A$ has to be set higher than the one resulting from physical measurements in order to ensure stable numerical solutions. On the other hand, the implementation of the scheme is very easy. Moreover, each of the employed numerical methods is very well known, meaning that possible issues and restrictions concerning the numerical schemes are already fully studied. Therefore, this full discretization provides reliable numerical solutions that can be used as benchmark to test the correctness of other numerical methods. 

We remark that~\eqref{eq1a} could be solved more efficiently in spectral space via fast cosine transform. Indeed, the cosine transform automatically imposes homogeneous Neumann boundary conditions. Moreover, spectral methods allow us to choose a relative low number of grid points and obtain high spatial accuracy, provided smooth data. Another possibility would be to assume periodic boundary conditions on $\rho$ and use the discrete Fourier transform also in the transverse coordinate. This would not affect heavily the numerical solution because the analysed region of the variable $\rho$ is much smaller than its total domain. Therefore, effects due to boundary conditions are negligible. In this work, we used a finite difference Crank--Nicolson scheme for its simplicity of implementation and a better comparison with the numerical solutions provided in~\cite{blanc11}. 

In the next session, we propose a different approach to~\eqref{eq0} with the final goal to obtain an highly parallelisable scheme. Moreover, the presented numerical scheme mitigates the stiffness inherent in~\eqref{eq1a} and  \eqref{eq1d}, and it is able to reproduce solutions with sharp gradients maintaining a reasonable size grid for $(\rho,\theta)$.

\subsection{Exponential integrators}
\label{subsec:ei}
Similarly to the splitting approach, instead of~\eqref{eq0} we consider the evolution equation 
\begin{equation}
\label{eq2a}
\begin{split}
\partial_{\sigma}V(\sigma,\rho,\theta) = \frac{1}{4\pi} & \partial_{\theta}^{-1} \left(\partial^2_{\rho} V(\sigma,\rho,\theta)\right) + \frac{B}{2}\,\partial_{\theta}V^2(\sigma,\rho,\theta) + \\
& A\,\partial^2_{\theta} V(\sigma,\rho,\theta) + 2\pi U_{\parallel}(\sigma,\rho) \partial_{\theta} V(\sigma,\rho,\theta) - U_{\perp}(\sigma,\rho)\partial_{\rho} V(\sigma,\rho,\theta).
\end{split}
\end{equation}
Here we treat $\partial^{-1}_{\theta}$ in frequency space. In particular, the antiderivative  $\partial^{-1}_{\theta}$ corresponds to the multiplication with $-\mathrm{i}/k$ in frequency space, more details are given in 
 section~\ref{subsubsec:fdei}.  We distinguish two parts on the right-hand side of~\eqref{eq2a}: linear terms with constant coefficients given by
\begin{equation}
\label{eq2}
\frac{1}{4\pi}\partial^{-1}_{\theta} \left(\partial_{\rho}^2 V\right) + A\partial^2_{\theta}V
\end{equation}
and linear terms with variable coefficient together with the nonlinear term 
\begin{equation}
\label{eq3}
2\pi U_{\parallel} \partial_{\theta} V -U_{\perp}\partial_{\rho} V + \frac{B}{2}\partial_{\theta} V^2.
\end{equation}
This partition motivates the use of exponential integrators. The advantage of these integrators lies in the fact that they integrate the terms in~\eqref{eq2} exactly. Therefore, the stiffness given by second derivatives vanishes. Let us define a linear operator $\mathcal{L}$ and a nonlinear operator $b$ by setting
\begin{equation}
\mathcal{L}(V) := \frac{1}{4\pi} \partial^{-1}_{\theta}\left(\partial_{\rho}^2 V \right) + A\partial^2_{\theta}V, \quad
b(\sigma,V):= 2\pi U_{\parallel} \partial_{\theta} V -U_{\perp}\partial_{\rho} V + \frac{B}{2}\partial_{\theta} V^2.
\end{equation}
Then,~\eqref{eq0} is rewritten fir short as
\[
\partial_{\sigma}V = \mathcal{L}(V) + b(\sigma,V). 
\]
Notice that we make the $\sigma$-dependence in $b$ explicit since~\eqref{eq3} depends on $U_{\perp}$ and $U_{\parallel}$ that are $\sigma$-dependent. 
The exact solution is obtained by using the variation of constants formula and it reads
\begin{equation}
\label{eq4}
V(\sigma + \Delta \sigma,\rho,\theta) = \mathrm{e}^{\Delta \sigma \mathcal{L}}V(\sigma,\rho,\theta) + \int_{0}^{ \Delta\sigma} \mathrm{e}^{(\Delta\sigma-s)\mathcal{L}}b(\sigma+s,V(\sigma + s,\rho,\theta))\,\mathrm{d}s.
\end{equation}
Notice that in~\eqref{eq4} is given in an implicit form only and  the analytical solution is not available, in general. The basic idea of exponential integrators is  to obtain numerical solutions by approximating the integral in~\eqref{eq4} with the available information, see~\cite{ostermann10}. The simplest (reasonable) approximation of the integral is done by replacing $b\big(\sigma+s,V(\sigma + s,\rho,\theta)\big)$ with $b\big(\sigma,V(\sigma,\rho,\theta)\big)$. Then, the $s$-dependence in the function $b$ is removed and we integrate $\mathrm{e}^{(\Delta\sigma-s)\mathcal{L}}$ exactly. This gives the exponential Euler method:
\[
V(\sigma + \Delta \sigma)\approx V^{n+1} = \mathrm{e}^{\Delta\sigma\mathcal{L}}V^{n} + \Delta\sigma\varphi_1(\Delta\sigma\mathcal{L})b(\sigma,V^n),
\] 
where $\varphi_1(z) = (\mathrm{e}^{z}-1)/z$ in an entire function. Notice that for brevity we suppressed the $(\rho,\theta)$ dependence of the variable $V$. The numerical scheme is of first order in $\sigma$, i.e. it has the same order of convergence of the splitting scheme described in section~\ref{subsec:spl}. Higher order exponential integrators can be constructed systematically.  We refer the reader to~\cite{ostermann10} for an exhaustive discussion about exponential integrators. In this work we use a two-stage second order exponential integrator (ExpRK22) similarly to~\cite{einkemmer20, einkemmer18}. The scheme is given by
\begin{equation}
\label{eq5}
\begin{split}
V^{n,*} &= \mathrm{e}^{\Delta\sigma\mathcal{L}}V^{n} + \Delta\sigma\varphi_1(\Delta\sigma\mathcal{L})b(\sigma,V^n), \\
V^{n+1} &= \mathrm{e}^{\Delta\sigma\mathcal{L}}V^{n} + \Delta\sigma\Big( \left(\varphi_1(\Delta\sigma\mathcal{L}) - \varphi_2(\Delta\sigma\mathcal{L}\right)b(\sigma,V^n) + \varphi_2(\Delta\sigma\mathcal{L})b(\sigma+\Delta\sigma,V^{n,*})\Big),
\end{split}
\end{equation}
where $\varphi_2(z) = (\mathrm{e}^z - 1 - z)/z^2$.  

\subsection{Full discretization with FFT and WENO}
\label{subsubsec:fdei}
In this section we describe the spatial discretisation which makes use of the fast Fourier transform (FFT) and the weighted essentially non-oscillatory (WENO) scheme. As usual, we consider a uniform discretization of the propagation distance 
\[
0\leq \sigma\leq\Sigma,\quad N_{\sigma}\in\N,\quad \Delta\sigma = \frac{\Sigma}{N_{\sigma}},\quad \sigma^n = n\Delta\sigma.
\]
The discretization in $(\rho,\theta)$ is done by pseudo-spectral methods. We adopt the same uniform discretization as in section~\ref{fdspl}. As it is well known, linear operators with constant coefficients can be treated efficiently in frequency space. On the other hand, variable coefficients and nonlinear terms should be computed in the physical space.  In contrast to section~\ref{subsec:spl} we assume periodic boundary conditions both in $\rho$ and $\theta$. Assuming periodic boundary conditions for $\rho$ does not affect the numerical simulations. Indeed, the domain is chosen large enough so that effects of the boundary conditions are negligible for the investigated region.

The numerical scheme~\eqref{eq5} is considered in frequency space to facilitate the computation of the operator $\mathrm{e}^{\Delta\sigma\mathcal{L}}$.  We apply the Fourier transform to both equations in~\eqref{eq5} and obtain
\begin{equation}
\label{eq11}
\begin{split}
\F(V^{n,*}) &= \mathrm{E}\odot\F(V^{n}) + \Delta\sigma\,\Phi_1\odot \F \left(b(\sigma^n,V^n)\right),\\
\F(V^{n+1}) &=  
\mathrm{E}\odot\F(V^{n})  + \Delta\sigma\Big( \left(\Phi_1 - \Phi_2\right)\odot\F\left(b(\sigma^n,V^n)\right) + \Phi_2\odot\F\left(b(\sigma^{n+1},V^{n,*})\right)\Big), 
\end{split}
\end{equation}
where 
\begin{align*}
\mathrm{E}_{jk} &= \exp(z_{jk}), \quad \Phi_{1,jk} = \frac{\exp(z_{jk}) - 1}{z_{jk}}, \quad
\Phi_{2,jk} =  \frac{\exp(z_{jk})-1-z_{jk}}{z^2_{jk}},
\end{align*}
\begin{align*}
z_{jk} = \Delta\sigma \left(-\frac{1}{4\pi}\frac{j^2}{\mathrm{i}k + \epsilon/4\pi} - A k ^2\right).
\end{align*}
The symbol $\odot$ denotes the component-wise product between two matrices, e.g. 
\[
\Big(\mathrm{E}\odot\F(V^{n})\Big)_{jk} = \mathrm{E}_{jk}\,\F(V^{n})_{jk}.
\]
The value $\epsilon$ in $z_{jk}$ corresponds to the machine epsilon, e.g. $\epsilon=2^{-53}$ for double precision floating point or $\epsilon=2^{-24}$ for single precision floating point. The definition of $z_{jk}$ is in accordance to the so called regularized Fourier multiplier. Roughly speaking, we can think of the regularized Fourier multiplier as a numerical trick in order to avoid treating the 0th $k$-frequency separately. A similar idea is used for example in~\cite{einkemmer20a,klein11}.

Notice that the computation of the terms involving $b$ is done in the physical space. This means that, starting from $\F(V^n)$ four discrete Fourier transforms are computed in order to obtain the numerical solution $\F(V^{n+1})$. Namely, one inverse and one forward Fourier transform to compute $\F(b(\sigma^n,V^n))$ in the first step and one inverse and forward Fourier transform to compute $\F(b(\sigma^{n+1},V^{n,*}))$ in the second step. 

What is left is the grid discretization of $b(\sigma^n,V^{n})$ and $b(\sigma^{n+1},V^{n,*})$ in the variables $(\rho,\theta)$, which is performed by the weighted essentially non-oscillatory scheme of order 5 (WENO5). This nonlinear numerical scheme has the advantage to limit oscillations in region where the solution is not regular, i.e., where the solution has sharp gradients, or it is even discontinuous. Moreover, in the region where the solution is smooth the WENO5 scheme reaches high order accuracy, in this case order 5. For an exhaustive introduction to ENO and WENO schemes we refer the reader to~\cite{shu98}. 

The WENO5 scheme is tailored to discretise gradients-like operators. Then, the solution is advanced in time by a chosen time integrator, in our case the ExpRK22 scheme.  
Notice that~\eqref{eq3} is not in a gradient-like form, indeed we have
\begin{equation}
\label{eq7}
\partial_{\theta} \left(-2\pi U_{\parallel} V - \frac{B}{2} V^2\right) + U_{\perp}\partial_{\rho} V = 0.
\end{equation}
In~\eqref{eq7} we used the fact that $U_{\parallel}$ is $\theta$-independent. However, we cannot write $U_{\perp}\partial_{\rho} V$ in gradient form because  $U_{\perp}=U_{\perp}(\sigma,\rho)$ is $\rho$-dependent. We rewrite the last term as $U_{\perp}\partial_{\rho} V = \partial_{\rho} (U_{\perp} V) - \partial_{\rho} U_{\perp}V$, so we obtain
\begin{equation}
\partial_{\sigma} V + \nabla\cdot f(V)  - \partial_{\rho} U_{\perp}V = 0,
\end{equation}
with 
\[
\quad \nabla = [\partial_{\rho},\partial_{\theta}],\quad f(V) = \left[-2\pi U_{\parallel} V - \frac{B}{2} V^2, U_{\perp} V\right]^T.
\]
Then, we discretize $\nabla\cdot f(V)$ by the WENO5 scheme. The extra term $-\partial_{\rho} U_{\perp}V$ does not require any approximation. The $\rho$-derivative of $U_{\perp}$ can be analytically computed (or numerically approximated) before starting the $\sigma$-evolution because the coefficient $U_{\perp}$ is known a priori for every $\sigma$ and $\rho$. Therefore, the extra term $-\partial_{\rho} U_{\perp}V$ reduces to the point-wise multiplication $-\partial_{\rho}U^{n}_{\perp,j}\,V_{jk}$, where $\partial_{\rho}U^{n}_{\perp,j}$ approximates$\partial_{\rho}U_{\perp}(\sigma^n,\rho_j)$ and $V_{jk}$ approximates  $V(\rho_j,\theta_k)$.

The WENO5 scheme in combination with the ExpRK22 scheme must fulfil a CFL condition induced by the terms collected in $b(\sigma, V)$ in order to provide stable solutions. We remark that $B$ and the variable coefficients $U_{\parallel}$, $U_{\perp}$ are relatively small ($\leq 5\cdot 10^{-2}$). This gives a mild CFL condition and we observe in numerical simulations a similar CFL condition as the one for the numerical scheme presented in section~\ref{subsec:spl}. 

The algorithm proceeds in two stages. In the first stage the quantity $\F(V^{n,*})$ is computed, given $\F(V^n)$, see~\eqref{eq11}. 
\begin{itemize}
\item The evaluation of $\mathrm{E}\odot \F(V^n)$ requires $\mathcal{O}(N_{\rho} N_{\theta})$ operations;
\item An inverse Fourier transform is performed via IFFT in order to get $V^n$ and compute $b(\sigma^n,V^n)$. IFFT works with $\mathcal{O}(N_{\rho}\log{N_{\rho}}N_{\theta}\log N_{\theta})$ operations and $b(\sigma^n,V^n)$ is computed via WENO5 in $\mathcal{O}(N_{\rho}N_{\theta})$ operations;
\item Then, a FFT is performed to obtain $\F(b(\sigma^n,V^n))$. 
\end{itemize}
This means that the first stage requires $\mathcal{O}(N_{\rho}\log{N_{\rho}}N_{\theta}\log N_{\theta})$ operations. 
Similar considerations apply to the second stage leading to the same asymptotic estimate.
This already gives an indication that the proposed numerical method outperforms (at least asymptotically) the one given in section~\ref{subsec:spl}. Moreover, WENO schemes together with FFT and IFFT are very suitable to parallelization. This results in a remarkable boost in terms of performance, as the numerical experiments in section~\ref{sec:ne} show. 

\section{Numerical results}
\label{sec:ne}
Numerical simulations are performed for the following parameters:
\[
0\leq\sigma \leq 120,\quad 0=\rho_{\min}\leq\rho\leq \rho_{\max}=400,\quad -13\pi=\theta_{\min} \leq \theta\leq 15\pi=\theta_{\max}.
\]
The intervals $[\rho_{\min},\rho_{\max}]$ and $[\theta_{\min},\theta_{\max}]$ are chosen large enough in order to mitigate the influence of the boundary conditions in the simulation. The variable $\rho$ is examined only in the interval $[133,267]$, while $\theta$ is examined in $[0,\theta_{\max}]$. A similar setting is also adopted in~\cite{blanc11}. 

\begin{figure}
\centering
\includegraphics[height=\textwidth, width=.9\textwidth]{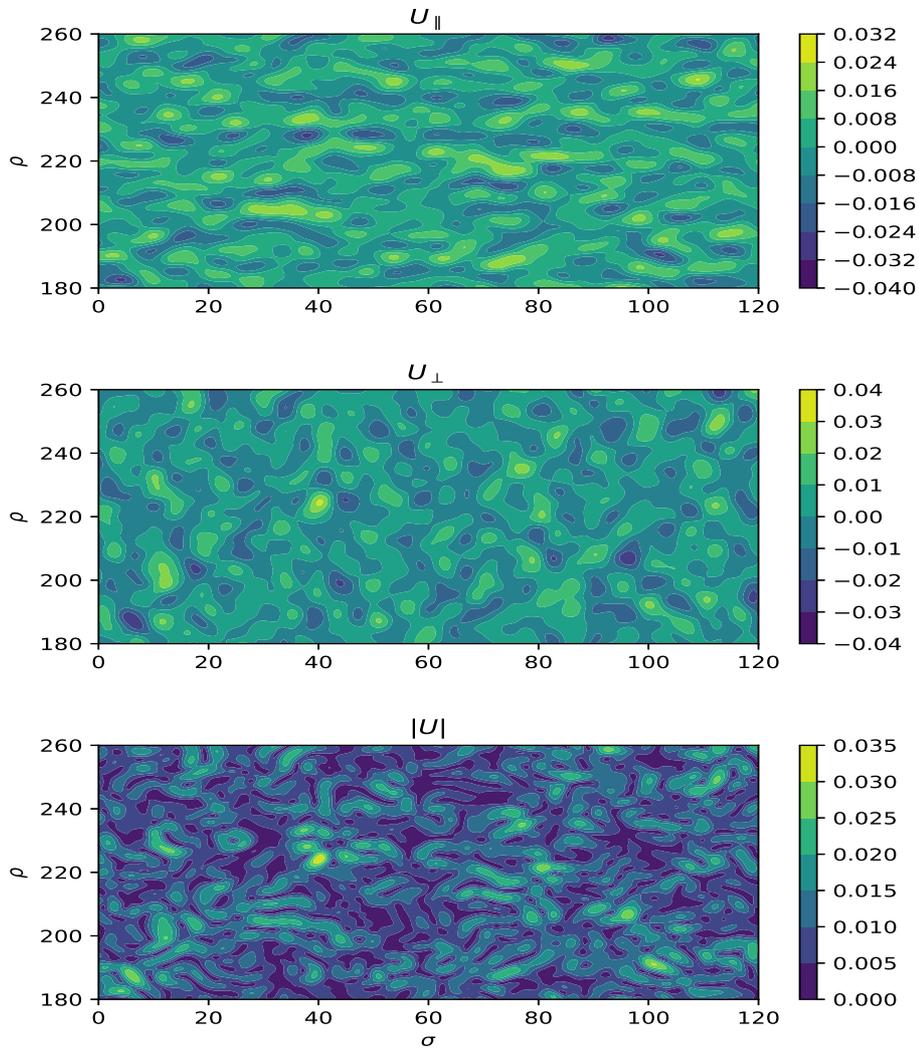}
\caption{On the $x$-axis the propagation distance $\sigma$. On the $y$-axis the transverse coordinate $\rho$. On top a density plot of first component of the velocity field $U$. In the middle a density plot of the second component of the velocity field $U$. Notice that the components are stretched along one direction. This is in accordance with physical observation, see~\cite{blanc11}. On the bottom the magnitude of $U$, notice that $\lvert U(\sigma,\rho)\rvert \leq 0.05$. }
\label{fig1}
\end{figure}

The variable coefficients $U_{\parallel}$, $U_{\perp}$ are generated before starting the simulation by using the procedure described in section~\ref{sec:mat}. An instance of  $U_{\parallel}$, $U_{\perp}$ is displayed in Fig~\ref{fig1}. For reasons of comparison we generate just one set of data $[U_{\parallel},U_{\perp}]^T$ for all numerical simulations. We use four different sets of values for $N_{\sigma}$, $N_{\rho}$ and $N_{\theta}$ collected in Table~\ref{tab0}. The inhomogeneous velocity fields are generated for $N_{\sigma} = 2400$ and $N_{\rho} = 10000$. For each of the remaining values in Table~\ref{tab0} the velocity fields $U_{\parallel}$, $U_{\perp}$ are sampled accordingly. 

\begin{figure}
\centering
\includegraphics[height=0.38\textwidth, width=0.78\textwidth]{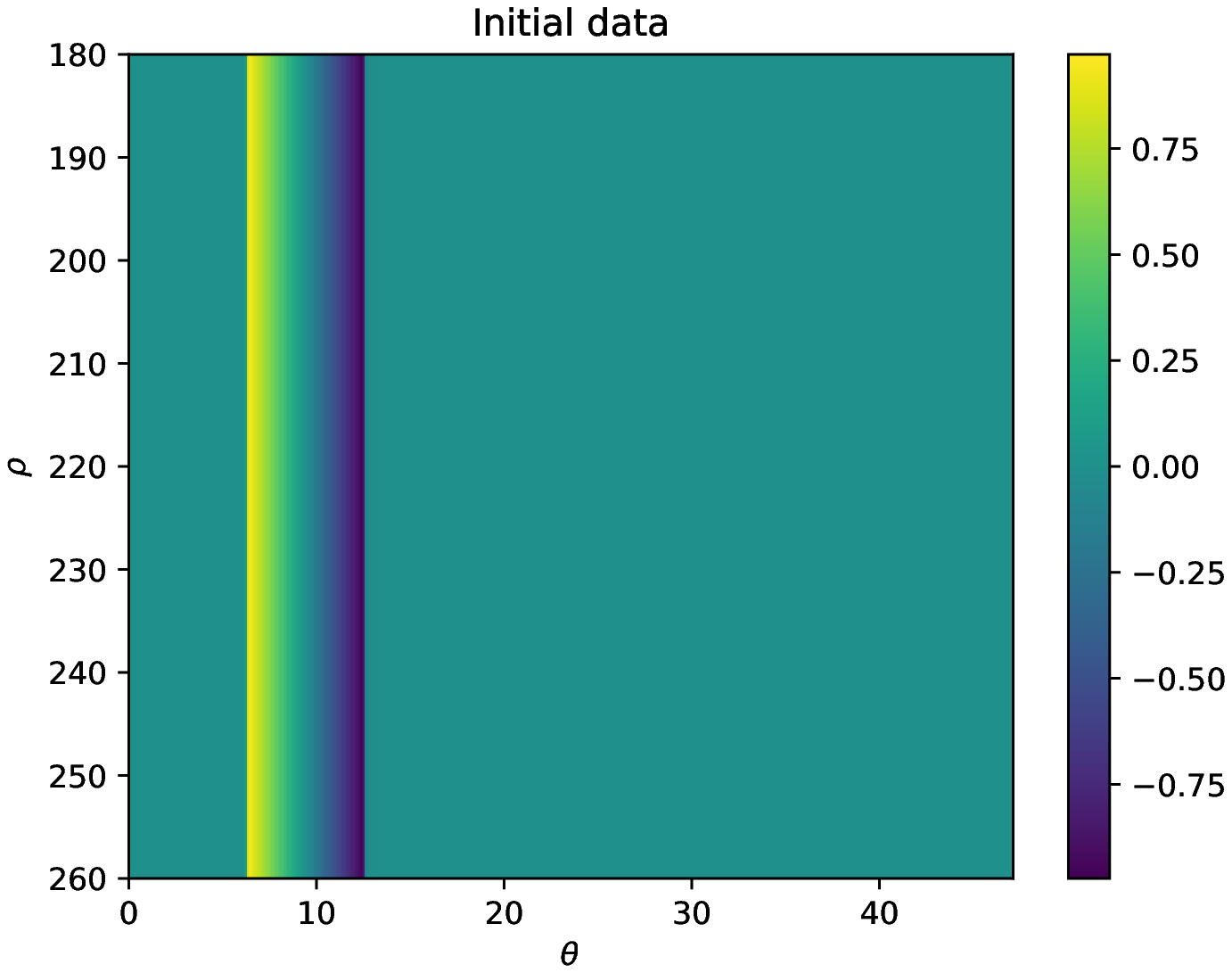}\\
\includegraphics[height=0.38\textwidth, width=0.78\textwidth]{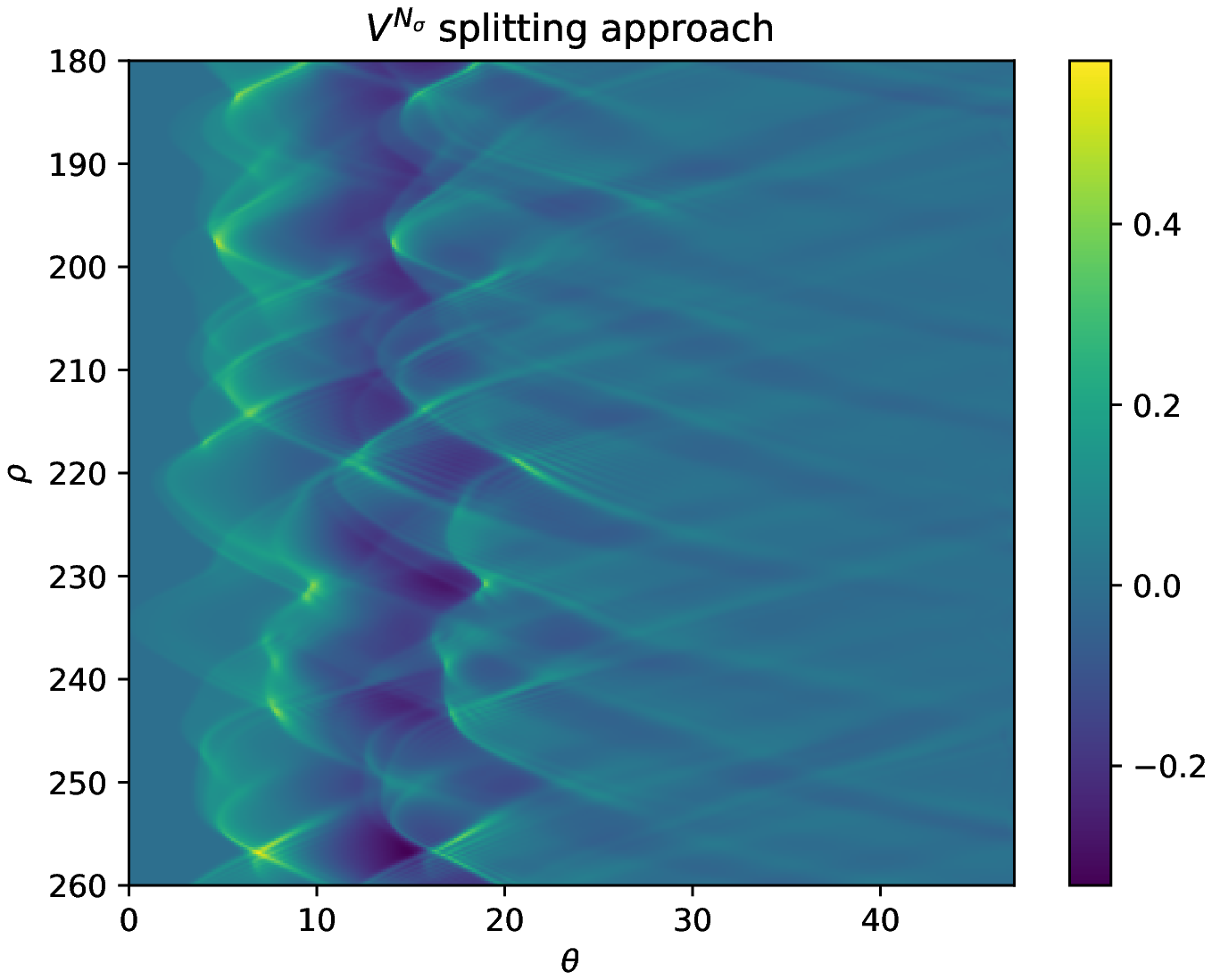} \\
\includegraphics[height=0.38\textwidth, width=0.78\textwidth]{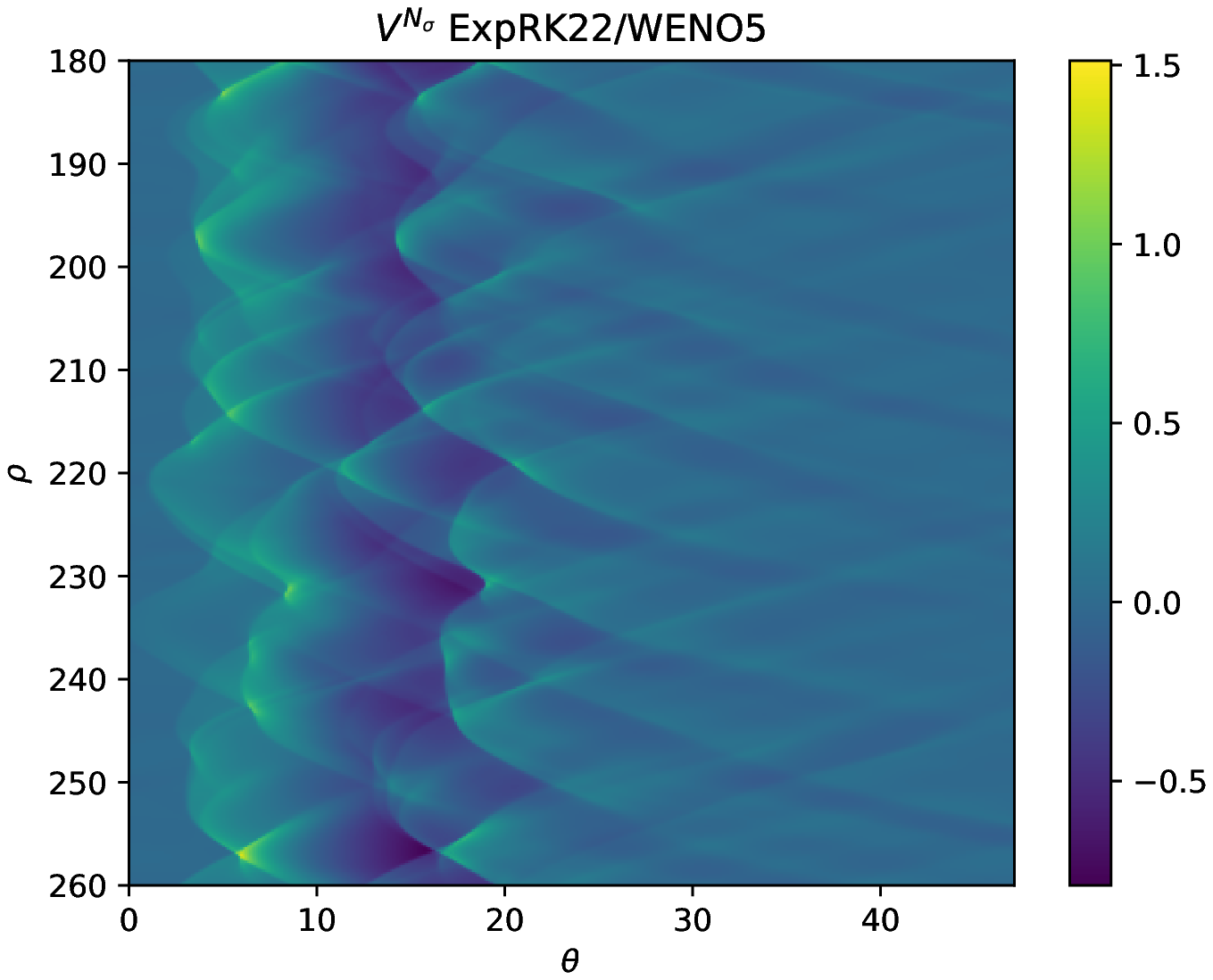}
\caption{Density plot of the initial value $V^0$ and the final solutions obtained by the splitting approach (middle) and the ExpRK22/WENO5 scheme (bottom). Notice that the final solutions have similar shape, but different amplitudes, the middle solution being three times smaller than the bottom solution.}
\label{fig0}
\end{figure}

The initial data is chosen as an $N$-wave pulse modeled by:
\[
V^0(\rho,\theta) = \frac{\theta-3\pi}{2\pi}\left(\tanh\left(\frac{B}{4A}(\theta-4\pi)\right) -\tanh\left({\frac{B}{4A}(\theta-2\pi)}\right)\right),\quad B=0.05.
\]
The absorption parameter $A$ will be discussed later in this section. In Fig.~\ref{fig0} we display the initial pulse $V^0$ and the final solution $V^{N_{\sigma}}$ obtained by the five way splitting scheme and the exponential integrator/WENO5 scheme. Here we employ $N_{\sigma} = 1200$, $N_{\rho} = 2500$ and $N_{\theta}=7\cdot 2^9$ for both schemes and display $V^{N_{\sigma}}$ for $(\rho,\theta)\in [133,267]\times [0,\theta_{\max}]$. The two solutions have a similar shape, but different amplitudes. The difference in amplitude is justified by the employed numerical methods. Indeed, the numerical schemes used in the splitting approach introduce numerical diffusion that is significantly higher than the numerical diffusion introduced by the WENO scheme. Therefore, the amplitude of the solution decreases much faster for the splitting approach than for the new approach presented in section~\ref{subsec:ei}.
The similar shape of the final solutions in Fig.~\ref{fig0}  confirms that assuming periodic boundary conditions on $\rho$ rather than homogeneous Neumann boundary conditions does not have appreciable effects on the final solution, at least not in the region of interest: $(\rho,\theta)\in[133,267]\times[0,\theta_{\max}]$.
 
\begin{table}[ht]
\centering
\begin{tabular}{c|ccc}
 & $N_{\sigma}$ & $N_{\rho}$ & $N_{\theta}$\Bstrut \\
\hline
Set 1 & $300$ & $1250$ & $7\cdot 64$\Tstrut \\
Set 2 & $600$ & $2500$ & $7\cdot 128$ \\
Set 3 & $1200$ & $5000$ & $7\cdot 256$ \\
Set 4 & $2400$ & $10000$ & $7\cdot 512$ \\
\end{tabular}
\caption{The four different sets of values for  $N_{\sigma}$, $N_{\rho}$ and $N_{\theta}$.}
\label{tab0}
\end{table}

\subsection{Splitting approach vs. exponential integrator approach}
We compare the two approaches by performing numerical simulations for the different values in Table~\ref{tab0}. These values are chosen so that the CFL conditions imposed by the splitting approach and the exponential integrator approach are fulfilled. In the splitting approach the CFL conditions imposed by the Godunov and the Lax-Wendroff scheme, used to solve~\eqref{eq1b},~\eqref{eq1e}, respectively, read
\[
N_{\theta}\leq \frac{28\pi}{30} N_{\sigma} \approx 3 N_{\sigma}\quad \text{and}\quad N_{\rho}\leq \frac{200}{3} N_{\sigma}\approx 66 N_{\sigma},
\]
where we used~\eqref{eq6} with the estimate $\lVert V\rVert_{\infty}\leq 5$. A similar CFL condition applies to the exponential integrator approach. We stress that these conditions are not restrictive, in particular the one imposed on $N_{\rho}$. This allows us to perform numerical simulations choosing reasonable grid sizes. 

\begin{figure}
\centering
\includegraphics[width=\textwidth , height=.9\textwidth]{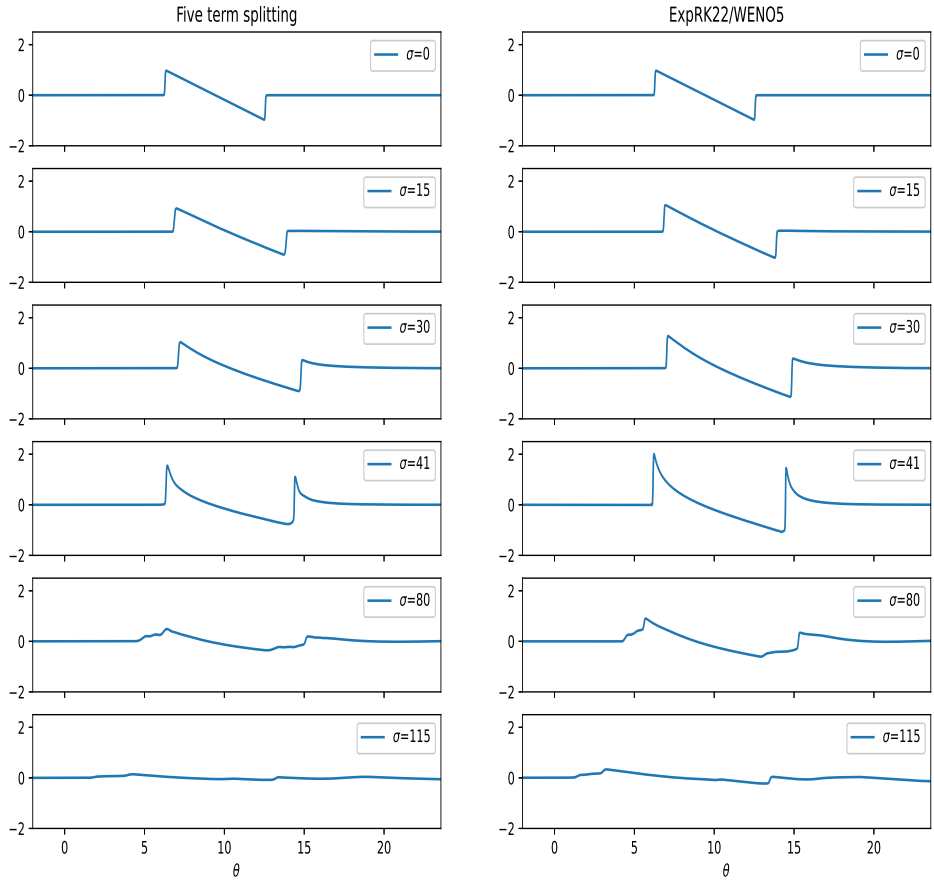}
\caption{Snapshots of the numerical solution $V$ at different propagation distances $\sigma$ along the transverse direction $\rho=260$. For $\sigma=30$ we notice the typical $U$-shape. As the propagation distance increases focusing effects are less likely and the original pulse gets dispersed.}
\label{fig3}
\end{figure}

In Fig.~\ref{fig3} we display the evolution of the numerical solutions obtained by the splitting approach and the exponential integrator approach at different propagation distances $\sigma$ along the fixed transverse coordinate $\rho =144$. Notice that for $\sigma=41$ the solutions exhibit a $U$-shape profile, similarly to~\cite{blanc11}. The profile of the solutions get flatter as $\sigma$ increases. This effect is due to the fact that the initial wave travels through the inhomogeneous medium following many different paths. This result in a scattering of the original pulse. Additionally, dissipative effects introduced by the numerical schemes in the splitting approach accentuate the flattening of the original pulse. 

Finally, in Table~\ref{tab4} we compare the total time needed to compute the numerical solution $V$ at the final propagation $\Sigma=120$ for the two approaches and for the different values in Table~\ref{tab0}. We observe that the splitting approach is faster for a small number of grid points, but is outperformed by a significant margin (up to a factor of 4 for large resolution)  by the ExpRK22 scheme as the grid gets more refined. This is in accordance with the theoretical estimates provided in section~\ref{sec:na}. Indeed, we have that the splitting scheme requires $\mathcal{O}( N_{\rho}N^2_{\theta})$ operations, while the exponential integrator requires $\mathcal{O}( N_{\rho}\log N_{\rho} N_{\theta}\log N_{\theta})$ operations. In addition, the exponential integrator/WENO scheme introduces less numerical diffusion and avoids oscillations as will be discussed in more detail below.

\begin{table}[ht]
\centering
\begin{tabular}{c|rr}
 & Splitting & ExpRK22 \Bstrut\\
\hline
Set 1 &  $35$ s &   $78$ s \Tstrut\\
Set 2 &  $569$ s & $649$ s \\
Set 3 & $11878$ s & $5196$ s \\
Set 4 & $>172800$ s & $42238$ s
\end{tabular}
\caption{Total time required for the simulations using different data set. The simulation times are reported in seconds. Data sets are reported in Table~\ref{tab0}. Notice that as the grid size gets refined the exponential integrator approach outperforms the splitting one. The maximum time (48 hours) was exceeded when the splitting approach was used for Set 4.}
\label{tab4}
\end{table}
\vspace{2mm}

\textbf{Absorption parameter}. The numerical solutions in the splitting approach are obtained by setting the absorption coefficient $A=3.4\cdot 10^{-4}$. This value is higher than the physical situation would warrant, where we have $A\sim 7\cdot 10^{-6}$. The reason of such a choice is  to prevent instabilities and/or oscillations of the numerical solutions, see the discussion in~\cite{blanc11}, section~\MakeUppercase{\romannumeral 2}-A. 

\begin{figure}
\centering
\includegraphics[width=\textwidth , height=.8\textwidth]{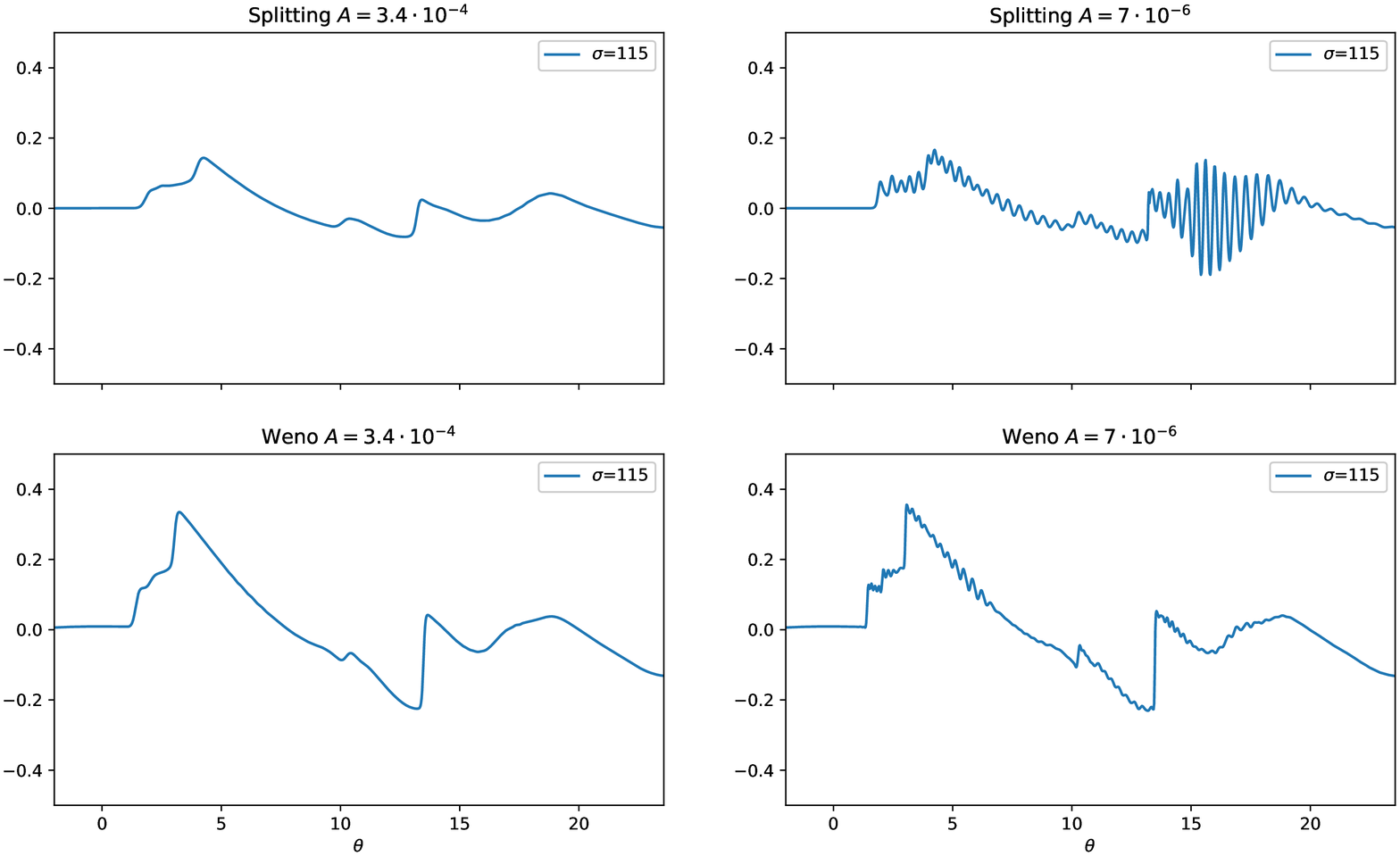}
\caption{The numerical solution $V$ as a function of $\theta$ at $(\sigma,\rho) = (115,144)$ for different values of the absorption parameter $A$. On top, the numerical solution obtained by the splitting approach, where oscillations are large for $A=7\cdot 10^{-6}$. At the bottom the solution obtained by the exponential integrator approach. We observe very small amplitude oscillations.}
\label{fig2}
\end{figure}

We test the new approach presented in section~\ref{subsec:ei} for both values of $A=3.4\cdot 10^{-4}$ and $A=7\cdot 10^{-6}$. In Fig.~\ref{fig2} we compare the numerical solution $V$ as a function of $\theta$ obtained by the two methods at $(\sigma,\rho)=(115,144)$. The discretisation is done by employing $N_{\sigma}=1200$, $N_{\rho} = 2500$ and $N_{\theta}=7\cdot 2^{9}$. Notice that oscillations are significantly bigger for the splitting approach, while they are negligible for the exponential integrator approach. The smaller oscillations in the exponential approach are due to use of the WENO scheme. Indeed, WENO schemes are able to capture shocks (i.e. where the solution has less regularity) by reducing the accuracy. However, in smooth regions the scheme recovers its precision, in the specific case the WENO scheme converges with order five. This allows us to choose the coefficient $A$ closer to the parameter given by the physical measurements.

\vspace{2mm}

\textbf{Convergence}. We test the convergence of the new approach, which is of second order in $\sigma$. To do so, we consider a reference solution $V_{\mathrm{ref}}$ obtained by using $N_{\sigma} = 2400$ points in the propagation direction and a final distance $\Sigma = 30$. Then, the reference solution is compared with numerical solutions $V_{\text{num}}$ obtained using a smaller number of points in the propagation direction. Both  $V_{\mathrm{ref}}$ and  $V_{\text{num}}$ are computed with the same number of grid points in the coordinate $\rho$ and $\theta$. In particular, we set $N_{\rho} = 5000$ and $N_{\theta} = 7\cdot 2^{8}$. 

We compute the relative error at the final propagation distance $\Sigma$ read off from
\[
\mathrm{err} = \frac{ \lVert V_{\text{ref}} - V_{\text{num}} \rVert_{{\ell}^2} }{\lVert V_{\text{ref}}\rVert_{\ell^2}},\quad 
\text{where}\quad  \lVert u \rVert_{\ell^2}:=\sqrt{\Delta\rho\Delta\theta \sum_{j,k} u_{j,k}^2}.
\] 
The convergence rate is given by the quotient of two consecutive error values. More specifically, let $N_1 < N_2$ be two different number of $\sigma$-points for the computation of two numerical solutions and $\mathrm{err}_1$,   $\mathrm{err}_2$ the associated relative errors. Then, the convergence rate $\beta$ is given by
\[
\frac{ \lVert \mathrm{err}_2 \rVert_{\ell^2} }{ \lVert \mathrm{err}_1\rVert_{\ell^2}} = \left(\frac{ N_2}{N_1}\right)^{-\beta}.
\]
In Table~\ref{tab1} we collect the relative errors with the convergence rate $\beta$. Notice that $\beta\approx 2$ indicates that the employed method is of second order. 

\begin{table}[ht]
\centering
\begin{tabular}{c|cccc|cc}
$N_{\sigma}$ & $\mathrm{err}$ & $\beta$ & & $N_{\sigma}$ & $\mathrm{err}$ & $\beta$\Bstrut\\
\cline{1-3}\cline{5-7}
$200$ & $1.728\mathrm{e}-03$ &  --  & & $600$ & $1.577\mathrm{e}-04$ & $2.20$\Tstrut \\
$300$ & $7.218\mathrm{e}-04$ & $2.15$ & & $800$ & $8.224\mathrm{e}-05$ & $2.26$\\
$400$ & $3.852\mathrm{e}-04$ &  $2.18$ & & $1200$ & $3.020\mathrm{e}-05$ & $2.47$\\

\end{tabular}
\caption{This table shows the relative errors and the convergence rate of the numerical solutions against a reference solution. The reference solution is computed by using $N_{\sigma}=2400$.}
\label{tab1}
\end{table}

\subsection{High performance computing}
We perform numerical simulations by using the exponential integrator method in combination with the WENO5 scheme presented in section~\ref{subsec:ei}.  

For the four different sets of values in Table~\ref{tab0} we measure the total time in seconds needed in order to compute the numerical solution at the final propagation distance $\Sigma~=~120$. The computer system used is an Intel Xeon scalable CPU Gold 6130 and a Titan V GPU. On the CPU the numerical simulations are performed by using 32 cores.

We carry out four tests. First, we test the sequential code. Second, we parallelize with the application programming interface OpenMP. Third, we parallelize with graphic processing units by using CUDA. In particular, for the GPU implementation we present two versions of the code: one in double precision floating point and one in single precision floating point. Performance results for the four different cases are reported in Table~\ref{tab2}. Numerical tests show a drastic speed up achieved using parallelisation with OpenMP compared to the sequential code, which is up to 22 times faster. The simulations are further accelerated when GPUs are involved. In particular, we observe that the single precision floating point implementations on GPUs run more than five times faster than the corresponding CPU implementations. The single precision implementation offers a speed up by a factor roughly two with respect to the double precision implementation, as expected.  To trade accuracy for precision might not always be a good choice. However, for this work single precision simulations still offer good results that give insight of the physical phenomena. We stress the fact that performance improvements of the simulations can be obtained only if the algorithm has a high rate of parallelisation. Differently, the simulation times on GPUs might result even inferior to the sequential code.     

An interesting aspect in HPC is to compare problems that are compute bound versus problems that are memory bound. Compute bound problems are of the kind that memory access is negligible with respect to the number of arithmetic operations, while the vice-versa holds for memory bound problems. This two aspects find place in the example treated in this work. More specifically the cost of the linear part are essentially memory bound, while the cost of the non-linear part are compute bound. We compare the time required to solve the linear parts against the WENO5 scheme. To illustrate this, let us consider the pseudo-code given in Algorithm~\ref{alg2}.

\begin{algorithm}
\caption{Exponential integrator/WENO5}
\begin{algorithmic}
\small
\State $V^0 = V^0(\rho_j,\theta_k),\quad 0\leq j < N_{\rho},\quad 0 \leq k < N_{\theta}$;
\For{$0\leq n < N_{\sigma}$}
\State STEP 1: Compute 
\begin{align*}
b(\sigma^n, V^n)\quad \text{with the WENO5 scheme;}
\end{align*} 
\State STEP 2: Compute 
\begin{align*}
\F(V^{n,*}) = \mathrm{E}\odot \F(V^n) + \Delta\sigma\, \Phi_1\odot \F(b(\sigma^n, V^n));
\end{align*}
\State STEP 3: Compute
\begin{align*}
b(\sigma^{n+1}, V^{n,*})\quad \text{with the WENO5 scheme;}
\end{align*}
\State STEP 4: Compute 
\begin{align*}
\F(V^{n+1}) = \mathrm{E}&\odot \F(V^n) \\
& + \Delta\sigma\, \{(\Phi_1-\Phi_2)\odot\F(b(\sigma^n, V^n)) + \Phi_2\odot\F(b(\sigma^{n+1}, V^{n,*}))\}.
\end{align*}
\EndFor
\end{algorithmic}
\label{alg2}
\end{algorithm}

STEP 1 and STEP 3 are responsible for the non-linear effects, while STEP 2 and STEP 4 for the linear ones. In Table~\ref{tab3} we report the required average time to compute the different steps in one iteration. To do so, we measure the total time to the completion of the simulation and divide it by $N_{\sigma}$. The results show how the computational cost of the linear effects is lower than the  non-linear ones for the sequential simulation. When the code is parallelised the situation is reversed. This is due to the fact that the WENO5 scheme is computationally bound, while the FFT is memory bound. 

\begin{table}[ht]
\centering
\begin{tabular}{c|rrrr}
Data set & \multicolumn{1}{c}{Seq.} & \multicolumn{1}{c}{OpenMP} & \multicolumn{1}{c}{CUDA double prec.} & \multicolumn{1}{c}{CUDA single prec.}\Bstrut\\
\hline
Set 1 &  $78.37$ s &  $3.776$ s (x$20.7$) &  $1.165$ s (x$3.2$) & $0.5564$ s (x$2.1$)\Tstrut\\
Set 2 &  $649$.4 s &  $29.02$ s (x$22.4$) &  $8.343$ s (x$3.5$) & $3.788$ s (x$2.2$)\\
Set 3 & $5196$ s & $232.6$ s (x$22.3$)  & $65.13$ s (x$3.6$) & $28.31$ s (x$2.3$)\\
Set 4 & $42238$ s &  $2775$ s (x$15.2$) &  $851.7$ s (x$3.3$) & $392.7$ s (x$2.2$)
\end{tabular}
\caption{Total simulation time in seconds for the different data sets in Table~\ref{tab0}. In brackets the speed-up factor with respect to the method employed on the previous column. E.g. for the Set 3  we have that the CUDA implementation in double precision is 3.6 times faster than the OpenMP implementation. The OpenMP implementation is 22.3 times faster than the Sequential version. This means that the CUDA implementation in single precision is almost 185 times faster than the sequential one.}
\label{tab2}
\end{table} 

\begin{table}[ht]
\centering
\begin{tabular}{c|cccc}
& Sequential  &   OpenMP & CUDA double prec. & CUDA single prec.\Bstrut  \\
\hline
nonlinear &  $9.834$ s &  $0.4319$ s &  $0.1654$ s &  $0.0778$ s\Tstrut\\
linear & $7.765$ s  & $0.7245$ s & $0.1894$ s & $0.0858$ s \\
\end{tabular}
\caption{Time required to simulate the non-linear (STEP 1+3) and linear (STEP 2+4) effects in a single iteration for the different code versions: sequential, Open MP, CUDA double precision and CUDA single precision. These results are obtained by using Set 4 in Table~\ref{tab0}. Notice how the computational cost of linear and non-linear effects scale equally with the parallelisation.}
\label{tab3}
\end{table}

\section{Conclusions}

This work is devoted to the study of the propagation of sonic-booms from the mid field into the far field. The mathematical model is given by a KZK-type equation which is a dispersive nonlinear partial differential equation. In the literature a numerical approach based on splitting methods has been proposed. In this work an adaptation of the algorithm is presented and discussed in details. One of the main disadvantages of the aforementioned approach consists in a difficult parallelisation of the algorithm. Therefore, we present a different approach based on exponential integrators in combination with WENO schemes. The new algorithm is highly parallelisable resulting in a tremendous acceleration (almost up to 185 times faster) with respect to its sequential version. Other than a reduced time of simulations, the exponential integrator approach brings additional benefits. The proposed algorithm achieves an higher accuracy with respect to the splitting approach, moreover the number of operations is (asymptotically) smaller. Finally, we observe numerically a significant reduction of oscillations (in amplitude) of the numerical solutions when the exponential integrator approach is used.  The new approach allows us to choose more grid-points and obtain numerical solutions that better describe the physical phenomena of $N$-wave propagation maintaining  the time of simulations relatively low. 

\clearpage

\bibliographystyle{siam}
\bibliography{references}

\end{document}